\documentclass[11pt,reqno]{amsproc}
\allowdisplaybreaks
\usepackage{cite}
\usepackage{tcolorbox}
\tcbuselibrary{breakable}
\usepackage{enumerate}
\usepackage{color}
\usepackage{mathrsfs}
\usepackage{fullpage}
\usepackage{graphicx}
\usepackage{subfigure}
\DeclareGraphicsExtensions{.eps}
\usepackage{ucs}
\usepackage[utf8x]{inputenc}
\usepackage{epstopdf}
\usepackage[semicolon,square,authoryear]{natbib}
\usepackage[debug=false, colorlinks=true, pdfstartview=FitV, 
linkcolor=blue, citecolor=blue, urlcolor=blue]{hyperref}

\usepackage{lineno}
\numberwithin{equation}{section}
\makeatletter
\def\maketag@@@#1{\hbox{\m@th\normalfont\normalsize#1}}
\makeatother

\usepackage{morefloats}
\newlength{\drop}
\definecolor{amethyst}{rgb}{0.6, 0.4, 0.8}
\definecolor{burgundy}{rgb}{0.5, 0.0, 0.13}

\title{\textbf{On boundedness and growth of unsteady solutions under the double porosity/permeability model}}

\author{\textbf{K.~B.~Nakshatrala} \\
  {\small Department of Civil and Environmental
    Engineering, University of Houston, Texas. \\
    \textbf{Correspondence to:}~knakshatrala@uh.edu}}


\begin{document}

\date{\today}

\begin{titlepage}
  \drop=0.1\textheight
  \centering
  \vspace*{\baselineskip}
  \rule{\textwidth}{1.6pt}\vspace*{-\baselineskip}\vspace*{2pt}
  \rule{\textwidth}{0.4pt}\\[\baselineskip]
       {\Large \textbf{\color{burgundy}
           On boundedness and growth of unsteady solutions under the double porosity/permeability model}}\\[0.3\baselineskip]
       \rule{\textwidth}{0.4pt}\vspace*{-\baselineskip}\vspace{3.2pt}
       \rule{\textwidth}{1.6pt}\\[\baselineskip]
       \scshape
       An e-print of the paper is available on arXiv:~1908.05771.  \par 
       \vspace*{\baselineskip}
       Authored by \\[\baselineskip]
       
       {\Large K.~B.~Nakshatrala\par}
       {\itshape Department of Civil \& Environmental Engineering \\
         University of Houston, Houston, Texas 77204--4003 \\ 
  \textbf{phone:} +1-713-743-4418, \textbf{e-mail:} knakshatrala@uh.edu \\
  \textbf{website:} http://www.cive.uh.edu/faculty/nakshatrala}\\[0.5\baselineskip]
       \vfill
       \begin{figure}[h]
       \centering
\includegraphics[scale=0.7]{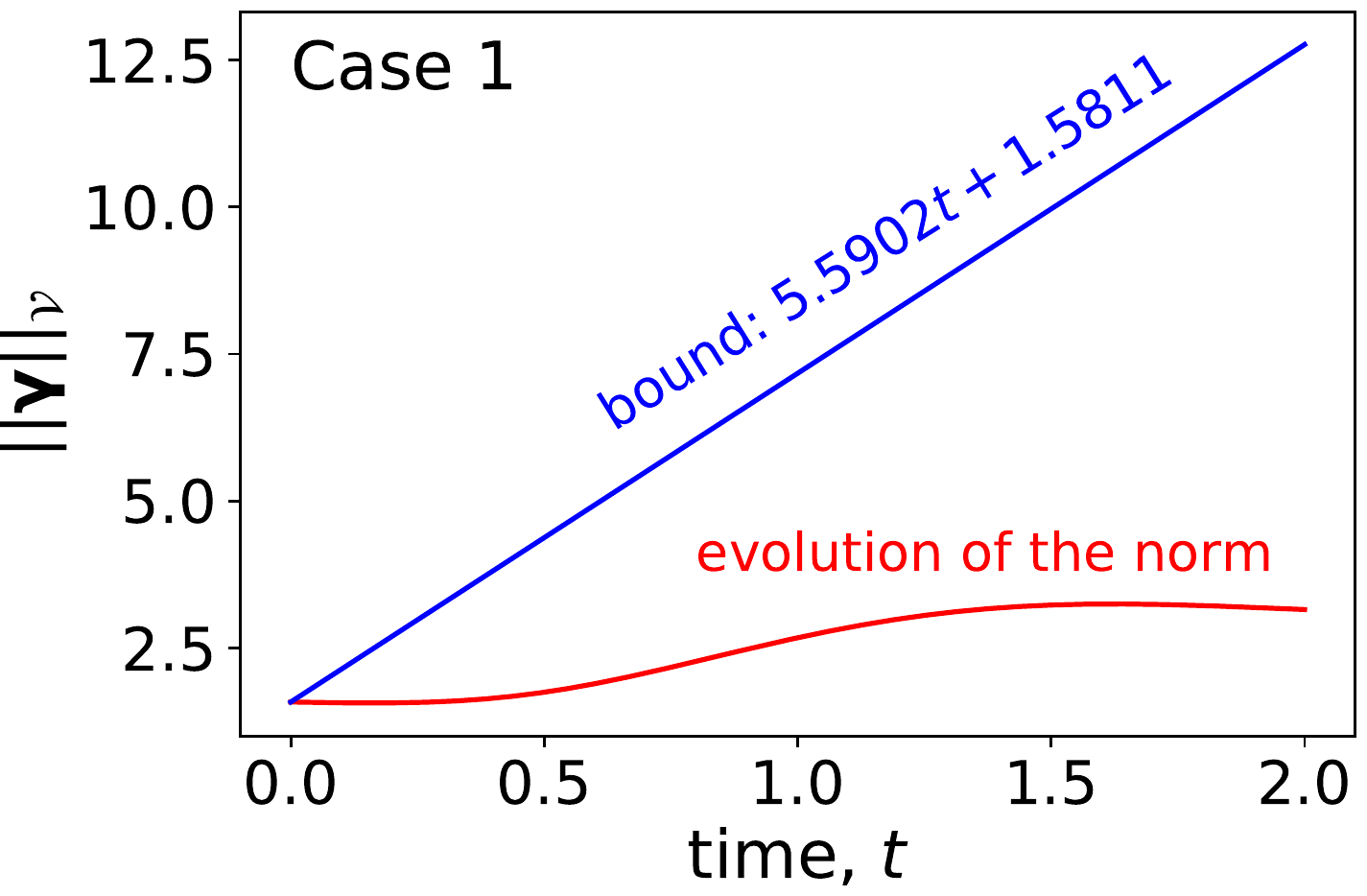}

\emph{This figure shows that the evolution of an unsteady solution under the DPP model satisfies the theoretical bound derived in this paper. $\|\boldsymbol{\Upsilon}\|_{\mathcal{V}}$ denotes a norm defined in terms of the velocities in the two pore-networks.}
\end{figure}
       \vspace*{\baselineskip}
           {\scshape 2020} \\
           {\small Computational \& Applied Mechanics Laboratory} \par
\end{titlepage}

\begin{abstract}
There is a recent surge in research activities on modeling the flow of fluids in porous media with complex pore-networks. A prominent mathematical model, which describes the flow of incompressible fluids in porous media with two dominant pore-networks allowing mass transfer across them, is the double porosity/permeability (DPP) model. However, we currently do not have a complete understanding of unsteady solutions under the DPP model. Also, because of the complex nature of the mathematical model, it is not possible to find analytical solutions, and one has to resort to numerical solutions. It is therefore desirable to have a procedure that can serve as a measure to assess the veracity of numerical solutions. In this paper, we establish that unsteady solutions under the transient DPP model are stable in the sense of Lyapunov. We also show that the unsteady solutions grow at most linear with time. These results not only have a theoretical value but also serve as valuable \emph{a posteriori} measures to verify numerical solutions in the transient setting and under anisotropic medium properties, as analytical solutions are scarce for these scenarios under the DPP model.
\end{abstract}

\maketitle

\section{OPENING STATEMENT}
\label{Sec:Bounded_Intro}
The study of the flow of fluids through porous media is central to various technological applications such as geological carbon sequestration, water purification, bioremediation, and enhanced oil recovery. Mathematical modeling often plays a crucial role in understanding the underlying dynamics in these applications, as the interior of the medium is inaccessible, notably, in subsurface applications. Also, the practical problems in these application areas are so complicated that they are not tractable via an analytical approach, and numerical solutions serve as the only viable tool.

Traditionally, Darcy equations have been used to model the flow of fluids through porous media. Because of the inherent simplicity of Darcy equations, and driven by their popularity, one can find in the literature analytical solutions for many problems \citep{strack2017analytical} and several robust numerical formulations \citep{chen2006computational,Masud_Hughes_CMAME_2002_v191_p4341}. Researchers have also established a myriad of mathematical properties that the solutions to Darcy equations satisfy, and these properties can serve as \emph{a posteriori} measures to assess the accuracy of numerical solutions; for example, see \citep{shabouei2016mechanics}. However, it is vital to realize that the Darcy model is valid under a plethora of assumptions \citep{Rajagopal_2007,nakshatrala2011numerical}. An assumption relevant to this paper is that the Darcy model assumes that the porous medium comprises a single dominant pore-network. 

Several emerging research endeavors have catalyzed a shift in attention towards porous media with complex pore-networks. The first one is the recent interest in exploring unconventional hydrocarbons, for example, oil and gas extraction from tight shale \citep{yao2013numerical}; shale exhibits double pore-networks with hydrocarbons trapped in low-permeability pore-network \citep{xu2013development}. Second, there is a growing interest in biomaterials to create improved medical devices, gain a deeper understanding of pathology, and facilitate biomimicry---creating new synthetic materials and systems that offer functionalities similar to living systems. Specifically, it is well-known that compact bone derives many of its mechanical and transport functionalities from its multiple pore-networks \citep{lemaire2006multiscale}. Last but not least, the latest revolution in manufacturing, such as 3D printing, has enabled us to build porous media with complicated pore structures and networks tailored to specific needs \citep{wang2016topological}. It should be clear that Darcy equations are not adequate to model the flow of fluids in porous media with such complex pore-networks. The flow dynamics in a porous medium with two pore-networks can be complex and differ from that of a single pore-network.

This inadequacy has resulted in the development of mathematical models, more complicated than Darcy equations, to address porous media with complex pore-networks. In particular, there has been a tremendous focus on flows in porous media with two or more dominant pore-networks with mass transfer across them; these works fall under the collective umbrella of double porosity/permeability (DPP) models. \citet{Barenblatt_Zheltov_Kochina_v24_p1286_1960_ZAMM} are often considered being the first to propose a DPP mathematical model. Through the years, there are several generalizations and different flavors of DPP models. The initial works considered steady-state flows; some notable ones include \citep{Barenblatt_Zheltov_Kochina_v24_p1286_1960_ZAMM,Warren_Root_1963_v3}. Subsequently, DPP model taking into account transient effects \citep{kazemi1969pressure,arbogast1989analysis} and deformation of the porous skeleton \citep{Khalili_v30_2003_Geophysical_research_letters,Borja_Koliji_2009,Choo_White_Borja_2015_IJG}. Analytical solutions to some simple transient problems are also obtained for DPP models, e.g., \citep{de1976analytic}. One can also find works that address derivation of these models using homogenization techniques \citep{Arbogast_Douglas_Hornung_1990,peszynska2009homogenization} and mathematical analysis (e.g., stability, existence of solutions) of these models \citep{arbogast1989analysis,hornung1990diffusion}.

However, the main assumption in the above-mentioned works is that the fluid is compressible with a constant coefficient of compressibility (e.g., see \citep[page 13]{arbogast1989analysis});  mathematically, $dp/d\rho = c/\rho$, where $p$ is the pressure, $\rho$ is the density of the fluid, and $c > 0$ is a constant. But in many applications, including the ones mentioned above, the density of the fluid does not change appreciably. Incompressibility of the fluid is an appropriate assumption in such situations. Notably, Darcy equations assume that the fluid is incompressible. Therefore, it is desirable to have a DPP model that generalizes Darcy equations, considers fluid's incompressibility, and applies to flow of fluids in porous media with double pore-networks. 

\citet{nakshatrala2018modeling} have developed such a DPP model, which will be the main focus of this paper; we will refer to this model as the DPP model from hereon. The mentioned paper also presents analytical solutions to the DPP model under steady-state conditions. However, analytical solutions for the DPP model under transient conditions and anisotropic medium properties are scarce. The primary reasons for the scarcity are (i) the mathematical model is complex, comprising four coupled partial differential equations expressed in terms of four field variables, and (ii) anisotropy gives rise to tensorial (rather than scalar) permeabilities, and tensorial quantities are more challenging to deal with than scalars. Other studies have developed numerical formulations to solve the governing equations under the DPP model \citep{joodat2018modeling, joshaghani2019stabilized}. Although some of these studies have addressed the transient DPP model, their focus has been narrow, primarily aimed at obtaining numerical solutions for specific initial-boundary value problems. To the best of the author's knowledge, there is no study on the general nature of unsteady solutions. It is fitting to recall the words of Truesdell in his book on Six Lectures on Modern Natural Philosophy \citep{truesdell1966method}: ``\emph{A mathematical theory is empty if it does not go beyond a few postulates, definitions, and routine calculations. Theorems must be proved, theorems, good theorems.}" Motivated by these words, this paper takes a modest step towards filling the lacuna in the theory of DPP. 

Returning to the other focus of this paper---regarding numerical solutions---it is imperative that a numerical simulator has to be well-tested by performing a series of checks before using it to carry out predictive numerical simulations. To put it another way, one needs to perform verification of solutions on the numerical simulator. Two popular strategies for verification are a comparison of the numerical solution with the analytical solution and the method of manufactured solutions \citep{Roache_2002_v124_p4_10_J_fluid_eng,WL_Oberkampf_AIAA_v36_p687}. However, as mentioned earlier, analytical solutions are scarce for the DPP model, and the method of manufactured solutions uses unrealistic boundary conditions and forcing functions. It is therefore desirable to have an alternate technique to check a numerical implementation so that one can use the formulation to solve other problems with confidence. It is also useful if we know the nature of the unsteady solutions and bounds on the growth or decay of solution fields with the time.

In the rest of this paper, we shall show that the solutions under the transient DPP model are stable in the sense of Lyapunov. We also show that the growth of the unsteady solutions can be at most linear in time under homogeneous boundary conditions. We will illustrate how one can utilize this mathematical result on the growth to construct a procedure to verify numerical solutions from a computer implementation.

\section{DPP MATHEMATICAL MODEL}
\label{Sec:S2_Unsteady_GE}
Let us consider a porous medium comprising two dominant pore-networks, referred to as the macro- and micro-pore networks. Each of these pore-networks has its hydromechanical properties; however, there could be a transfer of mass across the pore-networks. We denote the spatial domain by $\Omega \subset \mathbb{R}^{nd}$, where ``$nd$" denotes the number of spatial dimensions. A spatial point is denoted by $\mathbf{x}$. The divergence and gradient operators with respect to $\mathbf{x}$ are, respectively, denoted by $\mathrm{div}[\cdot]$ and $\mathrm{grad}[\cdot]$. We denote the time by $t \in [0,T]$, where $T$ denotes the length of the time interval of interest. 

For convenience, the quantities associated with the macro- and micro-pore networks will be indicated with subscripts 1 and 2, respectively. We denote the volume fractions by $\phi_1$ and $\phi_2$, the true (seepage) velocities by $\mathbf{v}_1(\mathbf{x},t)$ and $\mathbf{v}_2(\mathbf{x},t)$, the pressures by $p_1(\mathbf{x},t)$ and $p_2(\mathbf{x},t)$, and the bulk densities by $\rho_1$ and $\rho_2$. We denote the coefficient of viscosity and true density of the fluid by $\mu$ and $\gamma$, respectively. The bulk densities are related to the true density of the fluid as follows: 
\begin{align}
\label{Eqn:Unsteady_bulk_true_density}
\rho_1 = \gamma \phi_1 \quad \mathrm{and} \quad \rho_2 = \gamma \phi_2
\end{align}
It is also common to work in terms of the Darcy (discharge) velocities, which are defined as follows: 
\begin{align}
\mathbf{u}_1(\mathbf{x}) = \phi_1(\mathbf{x}) \mathbf{v}_{1}(\mathbf{x}) 
\quad \mathrm{and} \quad 
\mathbf{u}_2(\mathbf{x}) = \phi_2(\mathbf{x}) \mathbf{v}_{2}(\mathbf{x})
\end{align}
However, herein, we will work with the true velocities, and extending the framework based on the Darcy velocities is straightforward. 

The transient governing equations of the DPP model take the following form:
\begin{align}
  \label{Eqn:Unsteady_BoLM_1}
  &\rho_{1} \frac{\partial \mathbf{v}_1}{\partial t}
  + \mu \phi_1^{2} \mathbf{K}_{1}^{-1} \mathbf{v}_1
  + \phi_1 \mathrm{grad}[p_1] = \rho_1 \mathbf{b}_1 \\
  \label{Eqn:Unsteady_BoLM_2}
  &\rho_{2} \frac{\partial \mathbf{v}_2}{\partial t}
  + \mu \phi_2^{2} \mathbf{K}_{2}^{-1} \mathbf{v}_2
  + \phi_2 \mathrm{grad}[p_2] = \rho_2 \mathbf{b}_2 \\
  \label{Eqn:Unsteady_BoM_1}
  &\mathrm{div}[\phi_1 \mathbf{v}_1] = -\frac{\beta}{\mu} 
  (p_1 - p_2) \\
   \label{Eqn:Unsteady_BoM_2}
  &\mathrm{div}[\phi_2 \mathbf{v}_2] = +\frac{\beta}{\mu} 
  (p_1 - p_2)  
\end{align}
where $\mathbf{b}_1$ and $\mathbf{b}_2$ denote the specific body force in the pore-networks, and $\beta$ is a characteristic parameter of the porous medium. We often have $\mathbf{b}_1 = \mathbf{b}_2$ in practical situations; for example, the specific body force in each of the pore-networks is the acceleration due to gravity. It is important to note that $\phi_1$, $\phi_2$, $\rho_1$, $\rho_2$ and $\mu$ are all positive, and $\beta$ is non-negative. The permeabilities $\mathbf{K}_1$ and $\mathbf{K}_{2}$ are symmetric and positive definite tensors. The quantity $-\frac{\beta}{\mu}(p_1 - p_2)$ is the rate of volumetric transfer from the micro-pore network to the macro-pore network. The corresponding rate of mass transfer will then be $-\frac{\gamma \beta}{\mu}(p_1 - p_2)$.

The boundary conditions take the following form: 
\begin{subequations}
\begin{alignat}{2}
&\mathbf{v}_{1}(\mathbf{x},t) \cdot \widehat{\mathbf{n}}(\mathbf{x}) = v_{n1}(\mathbf{x},t) &&\quad 
\mathrm{on} \; \Gamma_{1}^{v} \\
&\mathbf{v}_{2}(\mathbf{x},t) \cdot \widehat{\mathbf{n}}(\mathbf{x}) = v_{n2}(\mathbf{x},t) &&\quad 
\mathrm{on} \; \Gamma_{2}^{v} \\
&p_1(\mathbf{x},t) = p_1^{\mathrm{p}}(\mathbf{x},t) &&\quad \mathrm{on} \; \Gamma_{1}^{p} \\
&p_2(\mathbf{x},t) = p_2^{\mathrm{p}}(\mathbf{x},t) &&\quad \mathrm{on} \; \Gamma_{2}^{p}
\end{alignat}
\end{subequations}
where $\Gamma_{1}^{v}$ and $\Gamma_1^{p}$ are the \emph{complementary} partitions of the boundary $\partial \Omega$, and likewise with $\Gamma_{2}^{v}$ and $\Gamma_2^{p}$ partitions. The initial conditions are prescribed as follows: 
\begin{align}
\mathbf{v}_{1}(\mathbf{x},0) = \mathbf{v}_1^{0}(\mathbf{x}) 
\quad \mathrm{and} \quad 
\mathbf{v}_{2}(\mathbf{x},0) = \mathbf{v}_2^{0}(\mathbf{x}) 
\quad \forall \mathbf{x} \in \Omega
\end{align}

We now mention two main assumptions behind the DPP model, and these assumptions are crucial in establishing the mathematical results presented in the subsequent sections. First, the volume fraction of each pore-network is independent of the time. That is, 
\begin{align}
\label{Eqn:Bounded_zeta_t_constant}
\frac{\partial \phi_1}{\partial t} = 0 
\quad \mathrm{and} \quad 
\frac{\partial \phi_2}{\partial t} = 0 
\end{align}
Second, the true density of the fluid, $\gamma$, is independent of the time. Note that this assumption is a stronger condition than the fluid is incompressible. In lieu of equations \eqref{Eqn:Unsteady_bulk_true_density} and \eqref{Eqn:Bounded_zeta_t_constant}, the bulk density in each pore-network is independent of time. That is, 
\begin{align}
\frac{\partial \rho_1}{\partial t} = 0 
\quad \mathrm{and} \quad 
\frac{\partial \rho_2}{\partial t} = 0 
\end{align}

A derivation along with a complete list of assumptions behind the DPP model are presented in \citep{nakshatrala2018modeling}.

\section{STABILITY}
We now show that the unsteady solutions under the transient DPP model are stable in the sense of a dynamical system. In particular, the solutions are Lyapunov stable \citep{dym2002stability,hale2012dynamics}. We assume the velocity boundary conditions to be homogeneous (i.e., $v_{n1} = 0$ on $\Gamma_{1}^{v}$ and $v_{n2} = 0$ on $\Gamma_{2}^{v}$). However, we allow the pressure boundary conditions to be non-homogeneous. 

\begin{figure}
\includegraphics[scale=0.65]{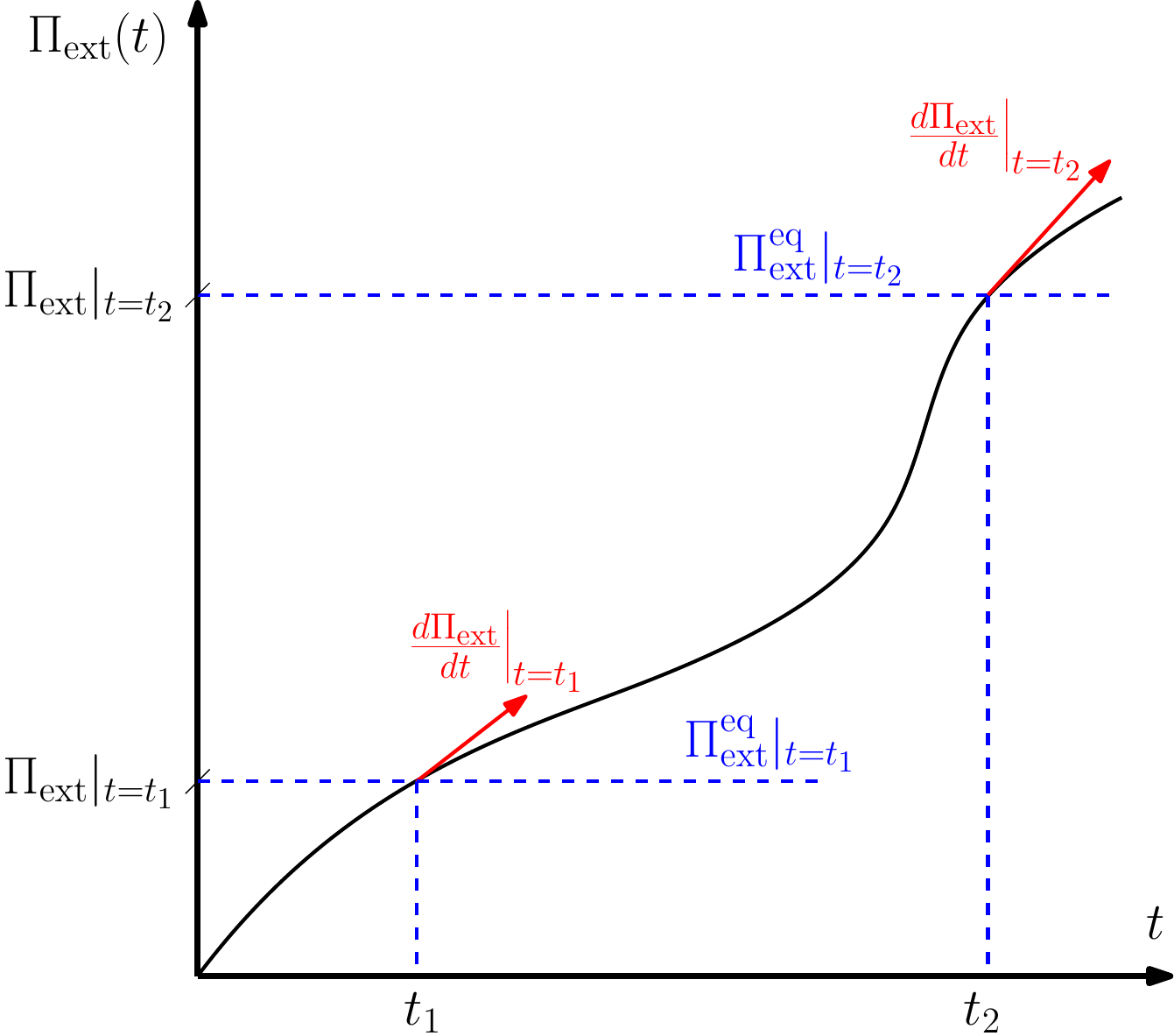}
\caption{This figure illustrates properties of $\Pi_{\mathrm{ext}}$ and $\Pi_{\mathrm{ext}}^{\mathrm{eq}}$, given by equations \eqref{Eqn:Unsteady_Pi_ext} and \eqref{Eqn:Unsteady_Pi_ext_dt}. We have denoted the one-parameter family of equilibrium states by $\Pi_{\mathrm{ext}}^{\mathrm{eq}}$.\label{Fig:Unsteady_Lyapunov}}
\end{figure}

For convenience, let 
\begin{align}
  \boldsymbol{\Upsilon} = \left\{\begin{array}{c}
  \mathbf{v}_1(\mathbf{x},t) \\
  \mathbf{v}_2(\mathbf{x},t)
  \end{array} \right\} 
\end{align}
We denote the equilibrium solution as follows:
\begin{align}
  \boldsymbol{\Upsilon}_{\mathrm{eq}} = \left\{
  \begin{array}{c}
    \mathbf{0} \\
    \mathbf{0}
  \end{array}\right\}
\end{align}
We consider the following functional, as a potential candidate for Lyapunov functional:
\begin{align}
  \mathbb{V}(\boldsymbol{\Upsilon})
  := \int_{\Omega} \left(\frac{1}{2} \rho_1 \mathbf{v}_1 \cdot
  \mathbf{v}_1 
  + \frac{1}{2} \rho_2 \mathbf{v}_2 \cdot
  \mathbf{v}_2 \right) \mathrm{d} \Omega
  + \Pi_{\mathrm{ext}} - \Pi_{\mathrm{ext}}^{\mathrm{eq}}
\end{align}
where $\Pi_{\mathrm{ext}}$ denotes the potential energy due to external loadings and $\Pi_{\mathrm{ext}}^{\mathrm{eq}}$ is the potential energy due to external loadings under an equilibrium state at a given instance of time. We assume the external loadings are conservative, which allows us to write the following: 
\begin{align}
\label{Eqn:Unsteady_definition_dPi_dt}
 \frac{d \Pi_{\mathrm{ext}}}{dt} = -\int_{\Omega} \left(\rho_1 \mathbf{b}_1 \cdot \mathbf{v}_1 
 + \rho_2 \mathbf{b}_2 \cdot \mathbf{v}_2 \right) \mathrm{d} \Omega
 +\int_{\Gamma_{1}^{p}} \phi_1 p_1^{\mathrm{p}} \mathbf{v}_1 \cdot \widehat{\mathbf{n}} \; \mathrm{d} \Gamma
 +\int_{\Gamma_{2}^{p}} \phi_2 p_2^{\mathrm{p}} \mathbf{v}_2 \cdot \widehat{\mathbf{n}} \; \mathrm{d} \Gamma
\end{align}
It is important to note that at a given instance of time, say $t = t_0$, we have 
\begin{align}
\label{Eqn:Unsteady_Pi_ext}
&\Pi_{\mathrm{ext}} \vert_{t = t_0}  = \Pi_{\mathrm{ext}}^{\mathrm{eq}} \vert_{t = t_0} 
\end{align}
but 
\begin{align}
\label{Eqn:Unsteady_Pi_ext_dt}
\left.\frac{d\Pi_{\mathrm{ext}}}{dt}\right|_{t = t_0}  \neq 0 \quad \mathrm{and} \quad 
\left.\frac{d\Pi_{\mathrm{ext}}^{\mathrm{eq}}}{dt}\right|_{t = t_0} = 0 
\end{align}
Figure \ref{Fig:Unsteady_Lyapunov} gives a pictorial description of the above equations. 

We now show that $\mathbb{V}$ is a Lyapunov functional for the DPP model. The task at hand is to show that the functional $\mathbb{V}$ satisfies the following three properties: 
\begin{enumerate}[(i)]
\item $\mathbb{V}(\boldsymbol{\Upsilon} = \boldsymbol{\Upsilon}_{\mathrm{eq}}) = 0$, 
\item $\mathbb{V}(\boldsymbol{\Upsilon} \neq \boldsymbol{\Upsilon}_{\mathrm{eq}}) > 0$, and 
\item $d\mathbb{V}/dt < 0$ for all $\boldsymbol{\Upsilon} \neq \boldsymbol{\Upsilon}_{\mathrm{eq}}$.
\end{enumerate}
The first two conditions are direct consequences of equation \eqref{Eqn:Unsteady_Pi_ext} and the definition of $\boldsymbol{\Upsilon}_{\mathrm{eq}}$. To establish the third condition, we proceed as follows: 
\begin{align}
  \frac{d\mathbb{V}}{dt} &= \int_{\Omega}
  \left(\rho_{1} \mathbf{v}_1 \cdot
  \frac{\partial \mathbf{v}_1}{\partial t}
  + \rho_{2} \mathbf{v}_2 \cdot
  \frac{\partial \mathbf{v}_2}{\partial t}
  \right) \mathrm{d} \Omega
  + \frac{d\Pi_{\mathrm{ext}}}{dt}
\end{align}
Using the balance of linear momentum in each pore-network, given by equations \eqref{Eqn:Unsteady_BoLM_1} and \eqref{Eqn:Unsteady_BoLM_2}, we obtain the following; 
\begin{align}
\frac{d\mathbb{V}}{dt} &=
-\int_{\Omega} \mathbf{v}_1 \cdot \left(\mu \phi_1^2 \mathbf{K}_1^{-1} \mathbf{v}_1 + \phi_1 \mathrm{grad}[p_1] - \rho_1 \mathbf{b}_1\right) 
  \mathrm{d} \Omega \nonumber \\
  &\quad -\int_{\Omega} \mathbf{v}_2 \cdot \left(\mu \phi_2^2 \mathbf{K}_2^{-1} \mathbf{v}_2 + \phi_2 \mathrm{grad}[p_2] - \rho_2 \mathbf{b}_2\right)   \mathrm{d} \Omega 
  + \frac{d\Pi_{\mathrm{ext}}}{dt}
\end{align}
Using the Green's identity and equation \eqref{Eqn:Unsteady_definition_dPi_dt}, we get the following: 
\begin{align}
  \frac{d\mathbb{V}}{dt}
  &= -\int_{\Omega} \left(\mathbf{v}_1 \cdot \mu \phi_1^2 \mathbf{K}_1^{-1} \mathbf{v}_1 + \mathbf{v}_2 \cdot \mu \phi_2^2 \mathbf{K}_2^{-1} \mathbf{v}_2 \right)  \mathrm{d} \Omega 
 +\int_{\Omega} \mathrm{div}[\phi_1 \mathbf{v}_1] p_1 \mathrm{d} \Omega
 +\int_{\Omega} \mathrm{div}[\phi_1 \mathbf{v}_2] p_2 \mathrm{d} \Omega \nonumber \\
 &\quad -\int_{\partial \Omega} \phi_1 p_1 \mathbf{v}_1 \cdot \widehat{\mathbf{n}} \; \mathrm{d} \Gamma
 -\int_{\partial \Omega} \phi_2 p_2 \mathbf{v}_2 \cdot \widehat{\mathbf{n}} \; \mathrm{d} \Gamma 
 +\int_{\Gamma_{1}^{p}} \phi_1 p_1^{\mathrm{p}} \mathbf{v}_1 \cdot \widehat{\mathbf{n}} \; \mathrm{d} \Gamma
 +\int_{\Gamma_{2}^{p}} \phi_2 p_2^{\mathrm{p}} \mathbf{v}_2 \cdot \widehat{\mathbf{n}} \; \mathrm{d} \Gamma
\end{align}
Using the boundary conditions for the pressures, and invoking the assumption that the velocity boundary conditions are homogeneous, we obtain the following: 
\begin{align}
  \frac{d\mathbb{V}}{dt}
  &= -\int_{\Omega} \left(\mathbf{v}_1 \cdot \mu \phi_1^2 \mathbf{K}_1^{-1} \mathbf{v}_1 + \mathbf{v}_2 \cdot \mu \phi_2^2 \mathbf{K}_2^{-1} \mathbf{v}_2 \right)  \mathrm{d} \Omega 
 +\int_{\Omega} \mathrm{div}[\phi_1 \mathbf{v}_1] p_1 \mathrm{d} \Omega
 +\int_{\Omega} \mathrm{div}[\phi_1 \mathbf{v}_2] p_2 \mathrm{d} \Omega 
\end{align}
Using the balance of mass in each pore-network, given by equations \eqref{Eqn:Unsteady_BoM_1} and \eqref{Eqn:Unsteady_BoM_2}, we obtain the following: 
\begin{align}
  \frac{d\mathbb{V}}{dt}
  &= -\int_{\Omega} \left(\mathbf{v}_1 \cdot \mu \phi_1^2 \mathbf{K}_1^{-1} \mathbf{v}_1 + \mathbf{v}_2 \cdot \mu \phi_2^2 \mathbf{K}_2^{-1} \mathbf{v}_2 \right)  \mathrm{d} \Omega 
 -\int_{\Omega} \frac{\beta}{\mu} (p_1 - p_2)^2 \mathrm{d} \Omega 
\end{align}
Noting that $\beta \geq 0$, $\mu > 0$, $\mathbf{K}_1$ and $\mathbf{K}_2$ are positive definite, we conclude that 
 \begin{align}
  \frac{d\mathbb{V}}{dt} < 0 \qquad \forall \boldsymbol{\Upsilon} \neq \boldsymbol{\Upsilon}_{\mathrm{eq}} 
  \end{align}
This implies that $\mathbb{V}$ is a non-increasing functional along the flow field, and this establishes that it is a Lyapunov functional for the DPP model. From the theory of dynamical systems \citep{luo2012stability}, we conclude that the solutions under the DPP model are Lyapunov stable.

\section{GROWTH OF UNSTEADY SOLUTIONS}
The governing equations, presented in Section \ref{Sec:S2_Unsteady_GE}, can be compactly written as the following constrained evolution problem: 
\begin{align}
\frac{d \boldsymbol{\Upsilon}}{dt} &= 
\mathcal{L}[\boldsymbol{\Upsilon},p_1,p_2] 
+ \mathbf{f} \\
\label{Eqn:Unsteady_G_constraint}
\mathbf{0} &= \mathcal{G}[\boldsymbol{\Upsilon},p_1,p_2]
\end{align}
where the linear operators $\mathcal{L}
[\cdot]$ and $\mathcal{G}[\cdot]$, and 
$\mathbf{f}$ are, respectively, defined 
as follows: 
\begin{align}
\label{Eqn:Unsteady_L_operator}
\mathcal{L}[\boldsymbol{\Upsilon},p_1,p_2] &= 
\left\{\begin{array}{c}
\left(-\mu \phi_1^{2} \mathbf{K}_{1}^{-1} 
\mathbf{v}_1 - \phi_1 \mathrm{grad}[p_1]\right) / \rho_1 \\
\left(-\mu \phi_2^{2} \mathbf{K}_{2}^{-1} 
\mathbf{v}_2 - \phi_2 \mathrm{grad}[p_2] \right / \rho_2
\end{array} \right\} \\
\label{Eqn:Unsteady_G_operator}
\mathcal{G}[\boldsymbol{\Upsilon},p_1,p_2] &= 
\left\{\begin{array}{c}
\mathrm{div}[\phi_1 \mathbf{v}_1] 
+ \frac{\beta}{\mu}(p_1 - p_2) \\
\mathrm{div}[\phi_2 \mathbf{v}_2] 
- \frac{\beta}{\mu}(p_1 - p_2) 
\end{array} \right\} 
\end{align}
\begin{align}
\mathbf{f} = 
\left\{\begin{array}{c}
\mathbf{b}_{1} \\
\mathbf{b}_{2} 
\end{array} \right\}
\end{align}
We assume homogeneous boundary conditions are enforced on the entire boundary. 

We denote the standard $L_2$ inner-product for scalar and vector fields defined on $\Omega$ by $\langle\cdot;\cdot\rangle$. That is, for given scalar fields $a$ and $b$ and vector fields $\mathbf{a}$ and $\mathbf{b}$ we have
\begin{align}
  \langle a;b\rangle 
  = \int_{\Omega} ab\;\mathrm{d} \Omega
  \quad \mathrm{and} \quad 
  \langle \mathbf{a};\mathbf{b}\rangle 
  = \int_{\Omega} \mathbf{a} \cdot \mathbf{b} 
  \; \mathrm{d} \Omega
\end{align}
We consider the following product function space: 
\begin{align}
  \mathcal{V} = \left(L_{2}(\Omega)\right)^{nd} \times \left(L_{2}(\Omega)\right)^{nd}
\end{align}
A natural inner-product on $\mathcal{V}$ will be 
\begin{align}
\langle \boldsymbol{\Upsilon};\widetilde{\boldsymbol{\Upsilon}}\rangle = \int_{\Omega} \left(\mathbf{v}_1 \cdot \widetilde{\mathbf{v}}_1 + \mathbf{v}_2 \cdot \widetilde{\mathbf{v}}_2 \right) \mathrm{d} \Omega \quad \forall \boldsymbol{\Upsilon}, \widetilde{\boldsymbol{\Upsilon}} \in \mathcal{V}
\end{align}
where 
\begin{align}
  \widetilde{\boldsymbol{\Upsilon}} = \left\{\begin{array}{c}
  \widetilde{\mathbf{v}}_1(\mathbf{x},t) \\
  \widetilde{\mathbf{v}}_2(\mathbf{x},t)
  \end{array} \right\}
\end{align}
The norm corresponding to the inner-product $\langle\cdot;\cdot\rangle$ is defined as follows: 
\begin{align}
  \|\boldsymbol{\Upsilon}\| := \sqrt{\langle \boldsymbol{\Upsilon};\boldsymbol{\Upsilon}\rangle}
\end{align}
However, noting that $\rho_1 > 0$ and $\rho_2 > 0$, we choose the following convenient inner-product on the function space $\mathcal{V}$:
\begin{align}
  \langle \boldsymbol{\Upsilon};\widetilde{\boldsymbol{\Upsilon}}\rangle_{\mathcal{V}} 
  = \int_{\Omega} \left( \rho_{1} \mathbf{v}_{1} 
  \cdot \widetilde{\mathbf{v}}_1 + \rho_{2} \mathbf{v}_{2} 
  \cdot \widetilde{\mathbf{v}}_{2} \right) \mathrm{d} \Omega 
  \quad \forall \boldsymbol{\Upsilon}, \widetilde{\boldsymbol{\Upsilon}} \in \mathcal{V}
\end{align}
The associated norm on $\mathcal{V}$ is defined as follows: 
\begin{align}
  \|\boldsymbol{\Upsilon}\|_{\mathcal{V}} := 
  \sqrt{\langle \boldsymbol{\Upsilon};
    \boldsymbol{\Upsilon}\rangle_{\mathcal{V}}}
\end{align}
Noting that the bulk densities are bounded below and bounded above by finite positive constants, the 
norms $\|\cdot\|$ and $\|\cdot\|_{\mathcal{V}}$ are equivalent. To wit, if 
\begin{align}
0 < \rho_{\mathrm{min}} \leq \rho_1, \rho_2 \leq \rho_{\mathrm{max}} < \infty
\end{align}
then we have 
\begin{align}
  \rho_{\mathrm{min}} \|\boldsymbol{\Upsilon}\| \leq \|\boldsymbol{\Upsilon}\|_{\mathcal{V}} \leq 
 \rho_{\mathrm{max}}  \|\boldsymbol{\Upsilon}\|
\end{align}

We first show that operator $\mathcal{L}$ is dissipative on $\mathcal{V}$ which is used then to establish that the growth of the unsteady solutions is at most linear with time. 

To establish that the operator $\mathcal{L}$ is dissipative on $\mathcal{V}$, we need to show the following:
\begin{align}
\label{Eqn:Unsteady_dissipative_inequality}
\langle \boldsymbol{\Upsilon} ; \mathcal{L}\rangle_{\mathcal{V}} \leq 0 \quad \forall \boldsymbol{\Upsilon} \in \mathcal{V}
\end{align}
We proceed by substituting the definition of $\mathcal{L}$, equation \eqref{Eqn:Unsteady_L_operator}, into the left hand side of \eqref{Eqn:Unsteady_dissipative_inequality}: 
\begin{align}
\langle \boldsymbol{\Upsilon} ; \mathcal{L}\rangle_{\mathcal{V}} = -\langle\mathbf{v}_1;\mu \phi_1^2 \mathbf{K}_1^{-1} \mathbf{v}_1\rangle - \langle\mathbf{v}_2;\mu \phi_2^2 \mathbf{K}_2^{-1} \mathbf{v}_2\rangle
-\langle\mathbf{v}_1;\phi_1 \mathrm{grad}[p_1]\rangle
-\langle\mathbf{v}_2;\phi_2 \mathrm{grad}[p_2]\rangle
\end{align}
Noting that $\mathbf{K}_1$ and $\mathbf{K}_2$ are positive definite tensors, $\mu > 0$, $\phi_1 > 0$ and 
$\phi_2 > 0$, we conclude the following: 
\begin{align}
\langle \boldsymbol{\Upsilon} ; \mathcal{L}\rangle_{\mathcal{V}} \leq 
-\langle\mathbf{v}_1;\phi_1 \mathrm{grad}[p_1]\rangle
-\langle\mathbf{v}_2;\phi_2 \mathrm{grad}[p_2]\rangle
\end{align}
Invoking Green's identity and noting that the boundary conditions are homogeneous, we obtain the following: 
\begin{align}
\langle \boldsymbol{\Upsilon} ; \mathcal{L}\rangle_{\mathcal{V}} \leq 
\langle\mathrm{div}[\phi_1 \mathbf{v}_1];p_1\rangle
+\langle\mathrm{div}[\phi_2\mathbf{v}_2];p_2\rangle
\end{align}
Using the incompressibility constraints, given by equations \eqref{Eqn:Unsteady_G_constraint} and \eqref{Eqn:Unsteady_G_operator}, we obtain the following: 
\begin{align}
\langle \boldsymbol{\Upsilon} ; \mathcal{L}\rangle_{\mathcal{V}} \leq -
\left\langle\frac{\beta}{\mu} \left(p_1 - p_2\right);\left(p_1 - p_2\right)\right\rangle
\end{align}
Noting that $\beta \geq 0$ and $\mu > 0$, we obtain the desired result: $\langle \boldsymbol{\Upsilon} ; \mathcal{L}\rangle_{\mathcal{V}} \leq 0$. 

We now address the growth of the unsteady solutions. We proceed as follows: 
\begin{align}
\|\boldsymbol{\Upsilon}\|_{\mathcal{V}} \frac{d}{dt}
  \left(\|\boldsymbol{\Upsilon}\|_{\mathcal{V}}\right) 
= \frac{d}{dt} \left(\frac{1}{2} \|\boldsymbol{\Upsilon}
\|^{2}_{\mathcal{V}}\right) 
&= \frac{d}{dt} \left(\frac{1}{2} \langle\boldsymbol{\Upsilon};\boldsymbol{\Upsilon}
\rangle_{\mathcal{V}}\right) \nonumber \\
&= \langle\boldsymbol{\Upsilon};\partial 
\boldsymbol{\Upsilon}/\partial t
\rangle_{\mathcal{V}} \nonumber \\
&= \langle \boldsymbol{\Upsilon};\mathcal{L}
[\boldsymbol{\Upsilon},p_1,p_2]\rangle_{\mathcal{V}} 
+ \langle \boldsymbol{\Upsilon};\mathbf{f}
\rangle_{\mathcal{V}} 
\end{align}
Noting that $\mathcal{L}$ is dissipative, inequality \eqref{Eqn:Unsteady_dissipative_inequality}, we obtain the following: 
\begin{align}
\|\boldsymbol{\Upsilon}\|_{\mathcal{V}} \frac{d}{dt}
  \left(\|\boldsymbol{\Upsilon}\|_{\mathcal{V}}\right) 
&\leq  \langle \boldsymbol{\Upsilon};\mathbf{f} \rangle_{\mathcal{V}} 
\end{align}
By invoking the Cauchy-Schwarz inequality, we obtain the following:
\begin{align}
\|\boldsymbol{\Upsilon}\|_{\mathcal{V}} \frac{d}{dt}\left(\|\boldsymbol{\Upsilon}\|_{\mathcal{V}}\right) \leq \|\boldsymbol{\Upsilon}\|_{\mathcal{V}} 
\|\mathbf{f}\|_{\mathcal{V}}
\end{align}
For $\|\boldsymbol{\Upsilon}\|_{\mathcal{V}} \neq 0$, we conclude:
\begin{align}
\frac{d}{dt}
  \left(\|\boldsymbol{\Upsilon}\|_{\mathcal{V}}\right) 
\leq \|\mathbf{f}\|_{\mathcal{V}}
\end{align}
By integrating both sides with time, we establish: 
\begin{align}
\label{Eqn:Unsteady_theoretical_bound}
\|\boldsymbol{\Upsilon}\|_{\mathcal{V}} \leq t f_{\mathrm{max}} + c
\end{align}
where 
\begin{align}
f_{\mathrm{max}} &= \mathop{\mathrm{max}}_{t \in [0,T]} \|\mathbf{f}(\mathbf{x},t)\|_{\mathcal{V}} \\
c &= \|\boldsymbol{\Upsilon}(\mathbf{x},t=0)\|_{\mathcal{V}}
\end{align}
Since $f_{\mathrm{max}}$ and $c$ are constants and finite, we conclude that the unsteady solutions under the DPP model grow at most linear with time if the driving forcing functions are bounded. 

\section{A REPRESENTATIVE NUMERICAL RESULT}
We now show how one can use the bound on the growth of unsteady solutions to check the veracity of numerical solutions. We proceed by outlining an initial-boundary value problem (IBVP)  under the DPP model and get numerical solutions of the IBVP using stable numerical formulation and discretization techniques. The norm $\|\boldsymbol{\Upsilon}\|_{\mathcal{V}}$ will be calculated using the resulting numerical solutions and compared with the derived theoretical bound on the norm. 

The computational domain is a unit square: $\Omega = (0,1) \times (0,1)$. The boundary conditions are no flow on the entire boundary for both the pore-networks (i.e., $v_{n1}(\mathbf{x}) = 0$ and $v_{n2}(\mathbf{x}) = 0$). 
Since the boundary conditions are no flow on the entire boundary for both the pore-networks (i.e., homogeneous velocity boundary conditions), we prescribe pressure at a point in the domain for one of the pore-networks to ensure uniqueness of solutions. For further details on the uniqueness of solutions, see \citep{nakshatrala2018modeling,joodat2018modeling}.

The time interval of interest for the numerical study is $[0,2]$. The backward Euler, which is an unconditionally stable time-stepping scheme, is used with a time-step of $0.001$. Table \ref{Table:Unsteady_params} provides the parameters used in the simulation. We considered two cases, each of which has a different specific body force. We chose different anisotropic permeabilities for the macro- and micro-pore networks, as we want the flow dynamics to be a characteristic of the DPP model and different from that of the Darcy model. See \citep{nakshatrala2018modeling} for a discussion on the scenarios under which Darcy equations can capture the solutions of the DPP model. 

\begin{figure}
\centering
\includegraphics[scale=0.2,clip]{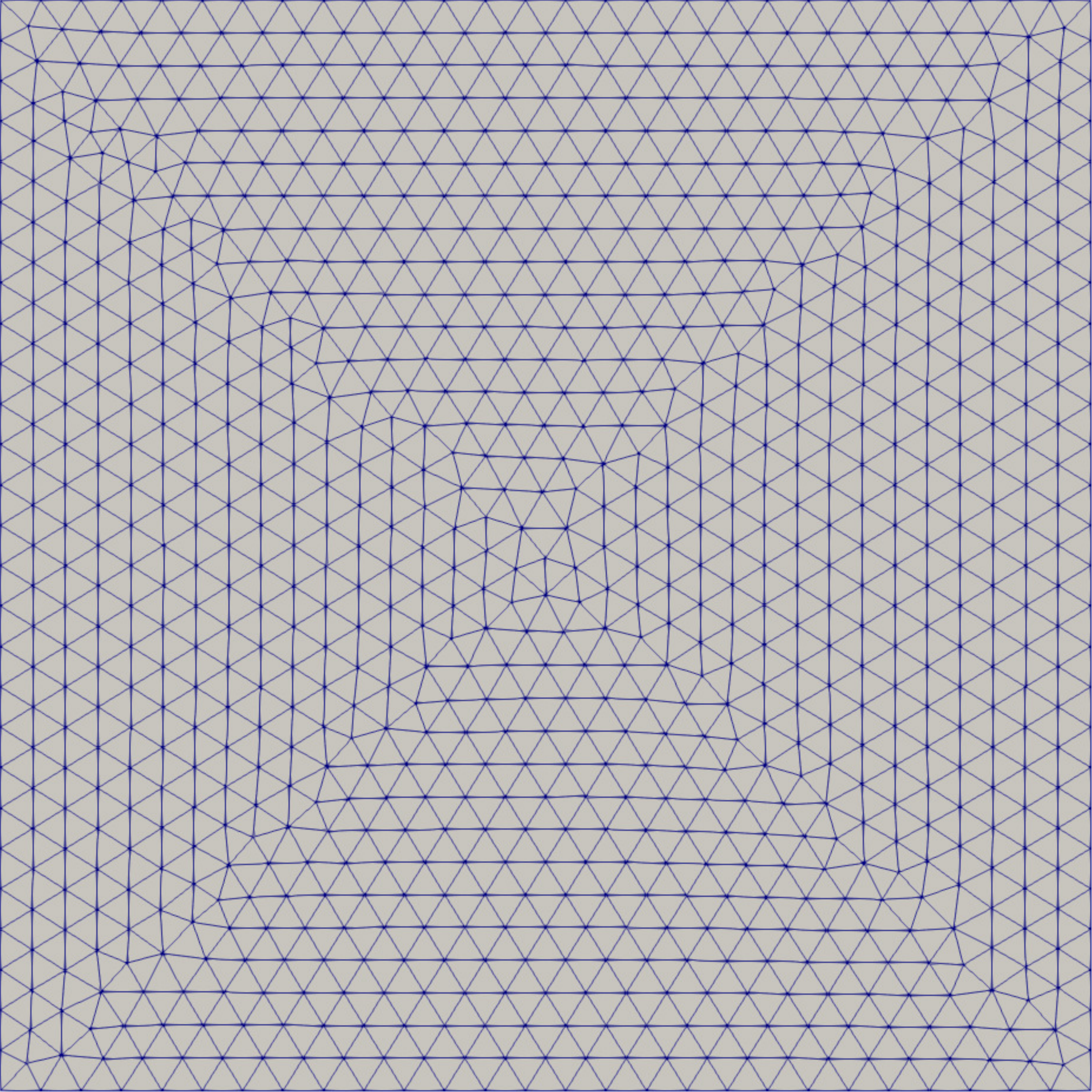}
\caption{This figure shows the mesh used in the numerical simulation.\label{Fig:Unsteady_mesh}}
\end{figure}

\begin{figure}
\centering
\subfigure[Case 1]{
\includegraphics[scale=0.17,clip]{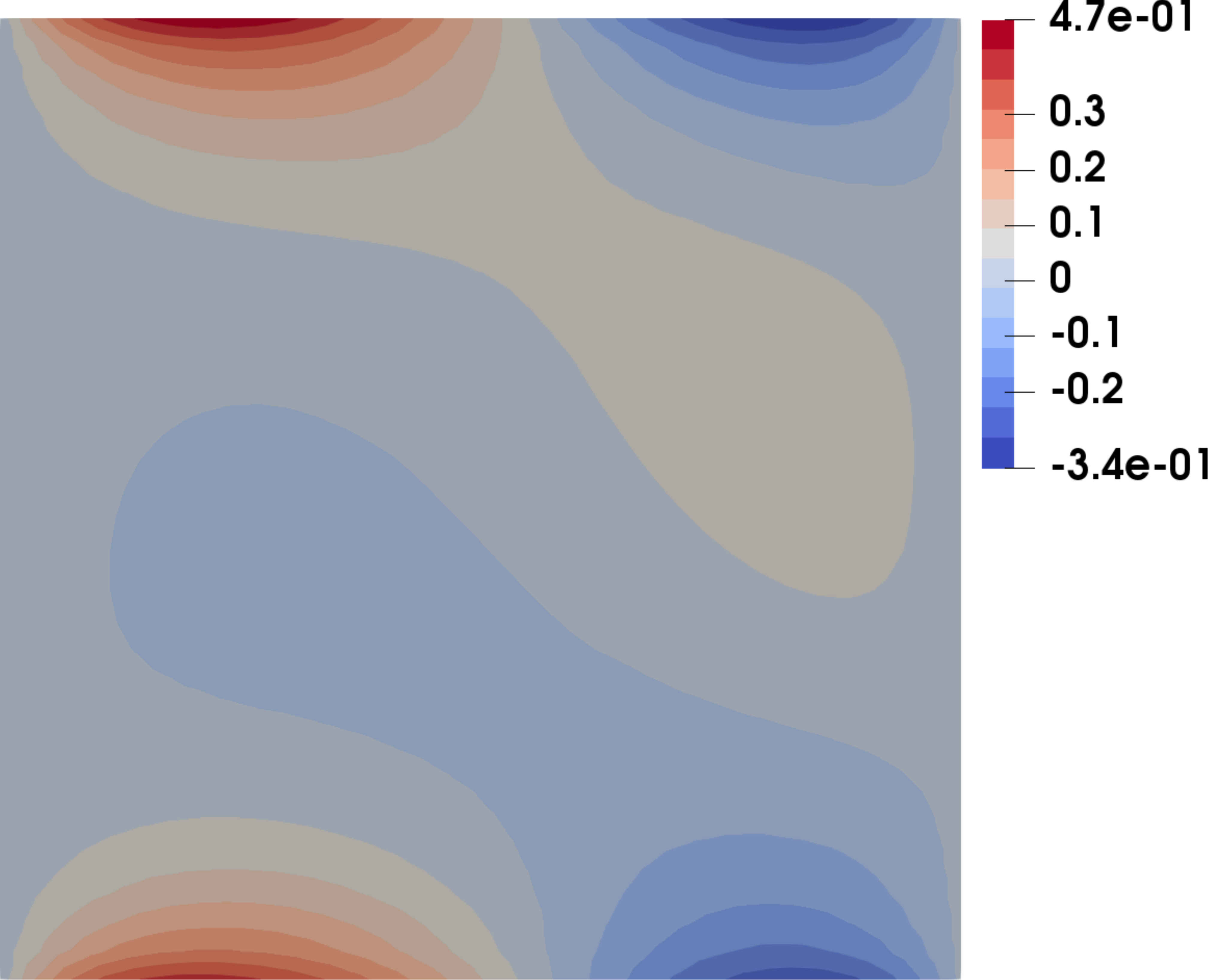}}
\subfigure[Case 2]{
\includegraphics[scale=0.17,clip]{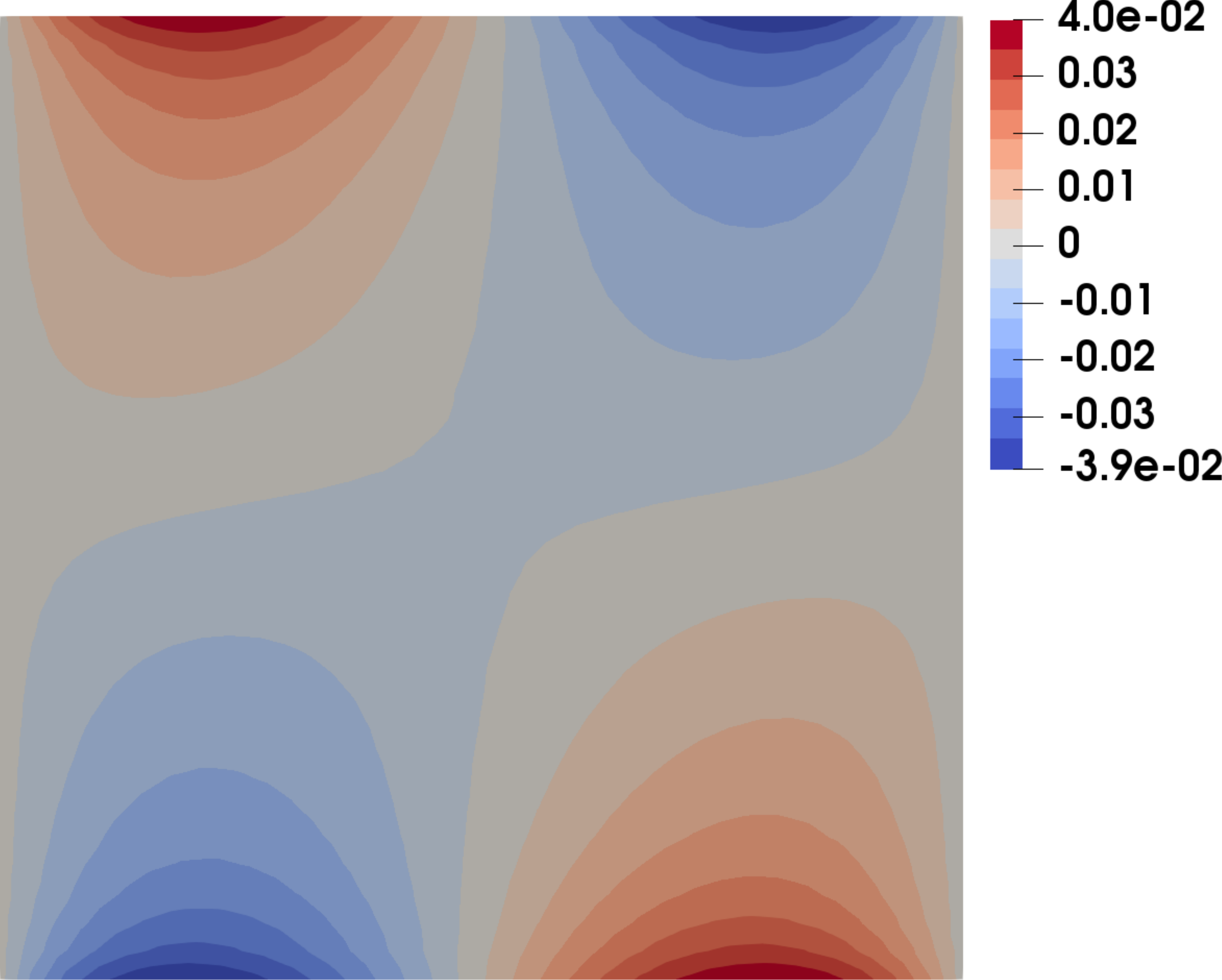}}
\caption{This figure shows the profiles of rates of volumetric transfer from the micro-pore network to the macro-pore network at time $t = 1$. The numerical results are stable and do not display any spurious oscillations which are typical of LBB violations. \label{Fig:Unstead_volumetric_transfer}}
\end{figure}

\begin{figure}
\subfigure{
\includegraphics[scale=0.53]{Figures/Evolution_of_norm_case_1.pdf}}
\subfigure{
\includegraphics[scale=0.53]{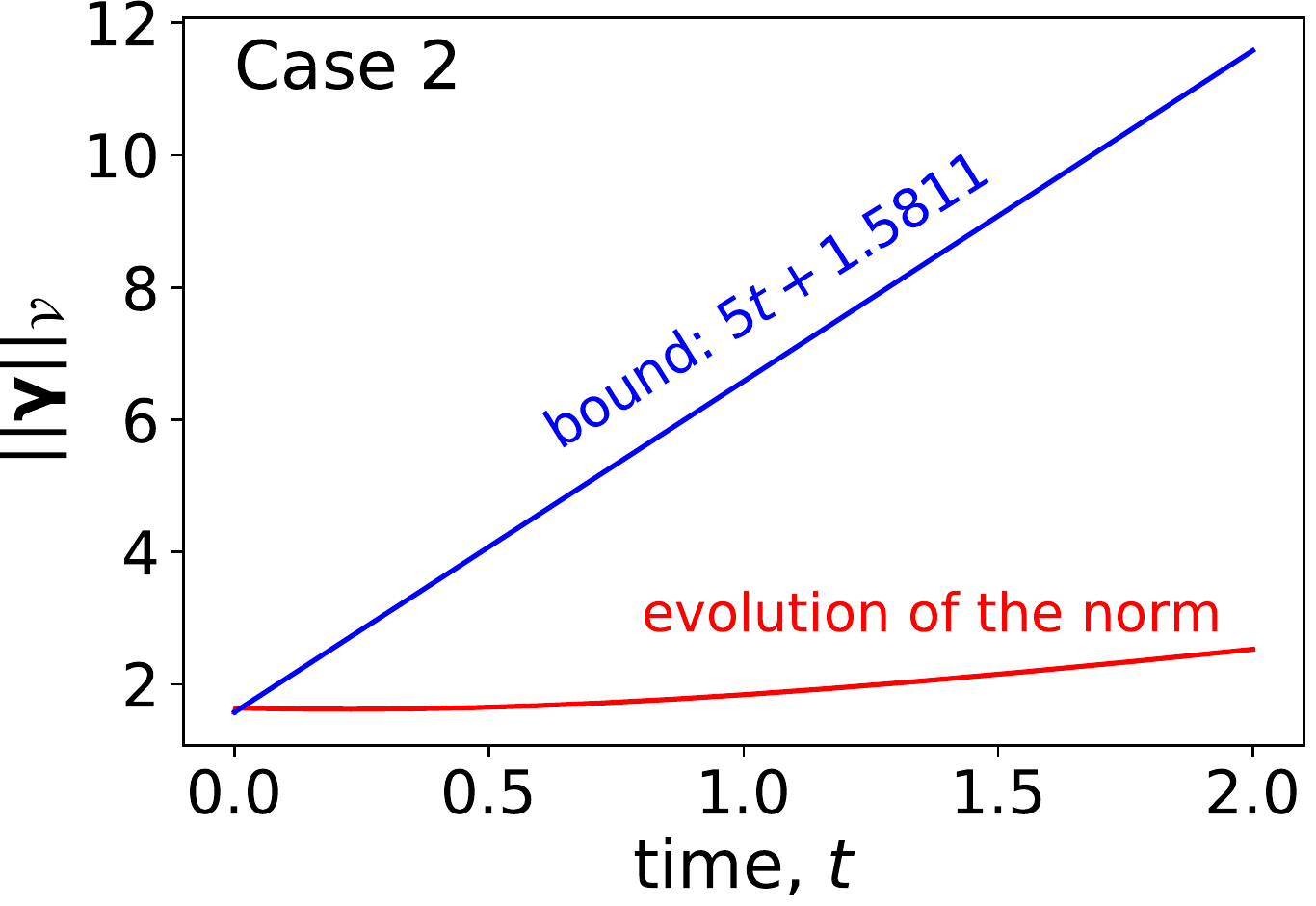}}
\caption{This figure shows the evolution of the norm $\|\boldsymbol{\Upsilon}\|_{\mathcal{V}}$ with time for the two cases considered in this paper. The corresponding bounds are plotted. The growth of the norm under the numerical simulations respect the theoretical bound derived in the paper.\label{Fig:Unstead_evolution_of_norm}}
\end{figure}

Figure \ref{Fig:Unsteady_mesh} shows the three-node triangular mesh used in the numerical simulation. We used the continuous Galerkin formulation with cubic interpolation for the velocity fields and linear interpolation for the pressure fields. This combination of interpolation functions---the so-called $P_3P_1$ interpolation on triangular elements---satisfies the Ladyzhenskaya-Babu\v ska-Brezzi (LBB) condition \citep{brezzi2012mixed}; see Figure \ref{Fig:Unstead_volumetric_transfer}. For the initial conditions, the intercept for the bound will be $c = 1.5811$. The slope of the bound for the two cases will be:
\begin{subequations}
\begin{align}
&\mbox{case 1:} \quad f_{\mathrm{max}} = 5.5902 \\
&\mbox{case 2:} \quad f_{\mathrm{max}} = 5
\end{align}
\end{subequations}
Figure \ref{Fig:Unstead_evolution_of_norm} shows that the numerical results satisfy the theoretical bound, given by equation \eqref{Eqn:Unsteady_theoretical_bound}.

\begin{table}[h]
\caption{Parameters used in the numerical simulations. \label{Table:Unsteady_params}}
\begin{tabular}{lc}\hline
\textbf{Quantity} & \textbf{Value} \\ \hline
True density, $\gamma$ & 1 \\
Coefficient of viscosity, $\mu$ & 1 \\
Mass transfer parameter, $\beta$ & 0.5 \\
Macro volume fraction, $\phi_1$ & 0.2 \\
Micro volume fraction, $\phi_2$ & 0.05 \\
Macro drag coefficient, $\mu \mathbf{K}_1^{-1}$ & $\left[\begin{array}{cc}
1 & 0.1 \\
0.1 & 0.9
\end{array}\right]$ \\
Micro drag coefficient, $\mu \mathbf{K}_2^{-1}$ & $\left[\begin{array}{cc}
100 & 5 \\
5 & 100
\end{array}\right]$ \\
Macro-velocity initial condition & $\mathbf{u}_1^{0} = \phi_1 \mathbf{v}_1^{0} = \left(\sin(\pi x) \cos(\pi y), -\cos(\pi x) \sin(\pi y)\right)$ \\
Micro-velocity initial condition & $\mathbf{u}_2^{0} = \phi_2 \mathbf{v}_2^{0} = \left(0, 0\right)$ \\
Case 1: specific body force, $\mathbf{b}_1=\mathbf{b}_2=\mathbf{b}$ & $\left(10\sin(\pi x t), 5\sin(2 \pi x y t)\right)$ \\
Case 2: specific body force, $\mathbf{b}_1=\mathbf{b}_2=\mathbf{b}$ & $\left(0, -10\right)$ \\ \hline
\end{tabular}
\end{table}

\section{CONCLUDING REMARKS}
We have presented two mathematical properties that the solutions under the transient DPP model satisfy: the unsteady solutions are Lyapunov stable, and they grow at most linear with time under homogeneous boundary conditions. Using a representative numerical example, we have shown that the second property---the nature of the growth---can serve as a verification procedure to check computer implementation of a numerical formulation. The attractive features are that the verification technique is easy to implement, in the form of \emph{a posteriori} measure, non-intrusive (i.e., one need not rewrite the computer code), and valid even under anisotropic medium properties. The results presented in this paper enlarge the repository of verification techniques to assess the accuracy of numerical simulations for the DPP model. 

One potential future work could be towards using the tools from functional analysis to get other mathematical properties and devise an array of verification techniques to assess the accuracy of numerical solutions under the transient DPP model.

\bibliographystyle{abbrvnat}
\bibliography{Master_References}

\begin{thebibliography}{31}
\providecommand{\natexlab}[1]{#1}
\providecommand{\url}[1]{\texttt{#1}}
\expandafter\ifx\csname urlstyle\endcsname\relax
  \providecommand{\doi}[1]{doi: #1}\else
  \providecommand{\doi}{doi: \begingroup \urlstyle{rm}\Url}\fi

\bibitem[Arbogast(1989)]{arbogast1989analysis}
T.~Arbogast.
\newblock Analysis of the simulation of single phase flow through a naturally
  fractured reservoir.
\newblock \emph{SIAM Journal on Numerical Analysis}, 26\penalty0 (1):\penalty0
  12--29, 1989.

\bibitem[Arbogast et~al.(1990)Arbogast, Douglas, and
  Hornung]{Arbogast_Douglas_Hornung_1990}
T.~Arbogast, J.~J. Douglas, and U.~Hornung.
\newblock {Derivation of the double porosity model of single phase flow via
  homogenization theory}.
\newblock \emph{SIAM Journal on Mathematical Analysis}, 21:\penalty0 823--836,
  1990.

\bibitem[Barenblatt et~al.(1960)Barenblatt, Zheltov, and
  Kochina]{Barenblatt_Zheltov_Kochina_v24_p1286_1960_ZAMM}
G.~I. Barenblatt, I.~P. Zheltov, and I.~N. Kochina.
\newblock {Basic concepts in the theory of seepage of homogeneous liquids in
  fissured rocks [strata]}.
\newblock \emph{Journal of Applied Mathematics and Mechanics}, 24:\penalty0
  1286--1303, 1960.

\bibitem[Borja and Koliji(2009)]{Borja_Koliji_2009}
R.~I. Borja and A.~Koliji.
\newblock {On the effective stress in unsaturated porous continua with double
  porosity}.
\newblock \emph{Journal of the Mechanics and Physics of Solids}, 57:\penalty0
  1182--1193, 2009.

\bibitem[Brezzi and Fortin(2012)]{brezzi2012mixed}
F.~Brezzi and M.~Fortin.
\newblock \emph{{Mixed and Hybrid Finite Element Methods}}, volume~15.
\newblock Springer Science \& Business Media, 2012.

\bibitem[Chen et~al.(2006)Chen, Huan, and Ma]{chen2006computational}
Z.~Chen, G.~Huan, and Y.~Ma.
\newblock \emph{Computational Methods for Multiphase Flows in Porous Media},
  volume~2.
\newblock SIAM Publishers, Philadelphia, 2006.

\bibitem[Choo et~al.(2016)Choo, White, and Borja]{Choo_White_Borja_2015_IJG}
J.~Choo, J.~White, and R.~I. Borja.
\newblock {Hydromechanical modeling of unsaturated flow in double porosity
  media}.
\newblock \emph{International Journal of Geomechanics}, 16\penalty0
  (6):\penalty0 D4016002, 2016.

\bibitem[de~Swaan(1976)]{de1976analytic}
A.~de~Swaan.
\newblock Analytic solutions for determining naturally fractured reservoir
  properties by well testing.
\newblock \emph{Society of Petroleum Engineers Journal}, 16\penalty0
  (03):\penalty0 117--122, 1976.

\bibitem[Dym(2002)]{dym2002stability}
C.~L. Dym.
\newblock \emph{{Stability Theory and its Applications to Structural
  Mechanics}}.
\newblock Dover Publications, Mineola, New York, 2002.

\bibitem[Hale and Ko{\c{c}}ak(2012)]{hale2012dynamics}
J.~K. Hale and H.~Ko{\c{c}}ak.
\newblock \emph{{Dynamics and Bifurcations}}, volume~3.
\newblock Springer Science \& Business Media, New York, 2012.

\bibitem[Hornung and Showalter(1990)]{hornung1990diffusion}
U.~Hornung and R.~E. Showalter.
\newblock Diffusion models for fractured media.
\newblock \emph{Journal of Mathematical Analysis and Applications},
  147\penalty0 (1):\penalty0 69--80, 1990.

\bibitem[Joodat et~al.(2018)Joodat, Nakshatrala, and
  Ballarini]{joodat2018modeling}
S.~H.~S. Joodat, K.~B. Nakshatrala, and R.~Ballarini.
\newblock {Modeling flow in porous media with double porosity/permeability: A
  stabilized mixed formulation, error analysis, and numerical solutions}.
\newblock \emph{Computer Methods in Applied Mechanics and Engineering},
  337:\penalty0 632--676, 2018.

\bibitem[Joshaghani et~al.(2019)Joshaghani, Joodat, and
  Nakshatrala]{joshaghani2019stabilized}
M.~S. Joshaghani, S.~H.~S. Joodat, and K.~B. Nakshatrala.
\newblock {A stabilized mixed discontinuous Galerkin formulation for double
  porosity/permeability model}.
\newblock \emph{Computer Methods in Applied Mechanics and Engineering},
  352:\penalty0 508--560, 2019.

\bibitem[Kazemi(1969)]{kazemi1969pressure}
H.~Kazemi.
\newblock Pressure transient analysis of naturally fractured reservoirs with
  uniform fracture distribution.
\newblock \emph{Society of Petroleum Engineers Journal}, 9\penalty0
  (04):\penalty0 451--462, 1969.

\bibitem[Khalili(2003)]{Khalili_v30_2003_Geophysical_research_letters}
N.~Khalili.
\newblock {Coupling effects in double porosity media with deformable matrix}.
\newblock \emph{Geophysical Research Letters}, 30:\penalty0 2153, 2003.

\bibitem[Lemaire et~al.(2006)Lemaire, Na{\"\i}li, and
  R{\'e}mond]{lemaire2006multiscale}
T.~Lemaire, S.~Na{\"\i}li, and A.~R{\'e}mond.
\newblock Multiscale analysis of the coupled effects governing the movement of
  interstitial fluid in cortical bone.
\newblock \emph{Biomechanics and Modeling in Mechanobiology}, 5\penalty0
  (1):\penalty0 39--52, 2006.

\bibitem[Luo et~al.(2012)Luo, Guo, and Morg{\"u}l]{luo2012stability}
Z.-H. Luo, B.-Z. Guo, and O.~Morg{\"u}l.
\newblock \emph{{Stability and Stabilization of Infinite Dimensional Systems
  with Applications}}.
\newblock Springer Science \& Business Media, London, 2012.

\bibitem[Masud and Hughes(2002)]{Masud_Hughes_CMAME_2002_v191_p4341}
A.~Masud and T.~J.~R. Hughes.
\newblock A stabilized mixed finite element method for {D}arcy flow.
\newblock \emph{Computer Methods in Applied Mechanics and Engineering},
  191:\penalty0 4341--4370, 2002.

\bibitem[Nakshatrala and Rajagopal(2011)]{nakshatrala2011numerical}
K.~B. Nakshatrala and K.~R. Rajagopal.
\newblock A numerical study of fluids with pressure-dependent viscosity flowing
  through a rigid porous medium.
\newblock \emph{International Journal for Numerical Methods in Fluids},
  67\penalty0 (3):\penalty0 342--368, 2011.

\bibitem[Nakshatrala et~al.(2018)Nakshatrala, Joodat, and
  Ballarini]{nakshatrala2018modeling}
K.~B. Nakshatrala, S.~H.~S. Joodat, and R.~Ballarini.
\newblock {Modeling flow in porous media with double porosity/permeability:
  Mathematical model, properties, and analytical solutions}.
\newblock \emph{Journal of Applied Mechanics}, 85\penalty0 (8):\penalty0
  081009, 2018.

\bibitem[Oberkampf and Blottner(1998)]{WL_Oberkampf_AIAA_v36_p687}
W.~L. Oberkampf and F.~G. Blottner.
\newblock {Issues in computational fluid dynamics: Code verification and
  validation}.
\newblock \emph{AIAA Journal}, 36:\penalty0 687--695, 1998.

\bibitem[Peszy{\'n}ska et~al.(2009)Peszy{\'n}ska, Showalter, and
  Yi]{peszynska2009homogenization}
M.~Peszy{\'n}ska, R.~Showalter, and S.-Y. Yi.
\newblock Homogenization of a pseudoparabolic system.
\newblock \emph{Applicable Analysis}, 88\penalty0 (9):\penalty0 1265--1282,
  2009.

\bibitem[Rajagopal(2007)]{Rajagopal_2007}
K.~R. Rajagopal.
\newblock {On a hierarchy of approximate models for flows of incompressible
  fluids through porous solids}.
\newblock \emph{Mathematical Models and Methods in Applied Sciences},
  17:\penalty0 215--252, 2007.

\bibitem[Roache(2002)]{Roache_2002_v124_p4_10_J_fluid_eng}
P.~J. Roache.
\newblock {Code verification by the method of manufactured solutions}.
\newblock \emph{Journal of Fluids Engineering}, 124:\penalty0 4--10, 2002.

\bibitem[Shabouei and Nakshatrala(2016)]{shabouei2016mechanics}
M.~Shabouei and K.~B. Nakshatrala.
\newblock Mechanics-based solution verification for porous media models.
\newblock \emph{Communications in Computational Physics}, 20\penalty0
  (5):\penalty0 1127--1162, 2016.

\bibitem[Strack(2017)]{strack2017analytical}
O.~D.~L. Strack.
\newblock \emph{{Analytical Groundwater Mechanics}}.
\newblock Cambridge University Press, Cambridge, U.K., 2017.

\bibitem[Truesdell(1966)]{truesdell1966method}
C.~A. Truesdell.
\newblock Method and taste in natural philosophy.
\newblock In \emph{Six Lectures on Modern Natural Philosophy}, pages 83--108.
  Springer, 1966.

\bibitem[Wang et~al.(2016)Wang, Xu, Zhou, Xu, Leary, Choong, Qian, Brandt, and
  Xie]{wang2016topological}
X.~Wang, S.~Xu, S.~Zhou, W.~Xu, M.~Leary, P.~Choong, M.~Qian, M.~Brandt, and
  Y.~M. Xie.
\newblock Topological design and additive manufacturing of porous metals for
  bone scaffolds and orthopaedic implants: A review.
\newblock \emph{Biomaterials}, 83:\penalty0 127--141, 2016.

\bibitem[Warren and Root(1963)]{Warren_Root_1963_v3}
J.~E. Warren and P.~J. Root.
\newblock {The behavior of naturally fractured reservoirs}.
\newblock \emph{Society of Petroleum Engineers Journal}, 3:\penalty0 245--255,
  1963.

\bibitem[Xu et~al.(2013)Xu, Haghighi, Li, and Cooke]{xu2013development}
B.~Xu, M.~Haghighi, X.~Li, and D.~Cooke.
\newblock Development of new type curves for production analysis in naturally
  fractured shale gas/tight gas reservoirs.
\newblock \emph{Journal of Petroleum Science and Engineering}, 105:\penalty0
  107--115, 2013.

\bibitem[Yao et~al.(2013)Yao, Sun, Fan, Wang, and Sun]{yao2013numerical}
J.~Yao, H.~Sun, D.-Y. Fan, C.-C. Wang, and Z.-X. Sun.
\newblock Numerical simulation of gas transport mechanisms in tight shale gas
  reservoirs.
\newblock \emph{Petroleum Science}, 10\penalty0 (4):\penalty0 528--537, 2013.

\end{thebibliography}
\end{document}